\documentclass{agtart_a}
\pdfoutput=1

\usepackage{graphicx}

\title{Heegaard splittings and the pants complex}

\author{Jesse Johnson}
\givenname{Jesse}
\surname{Johnson}
\address{Mathematics Department\\
University of California\\\newline
Davis, CA 95616\\USA}
\email{jjohnson@math.ucdavis.edu}
\urladdr{}

\volumenumber{6}
\issuenumber{}
\publicationyear{2006}
\papernumber{32}
\startpage{853}
\endpage{874}

\doi{}
\MR{}
\Zbl{}



\received{5 May 2006}
\revised{}
\accepted{11 May 2006}
\published{11 July 2006}
\publishedonline{11 July 2006}
\proposed{}
\seconded{}
\corresponding{}
\editor{}
\version{}

\arxivreference{math.GT/0509680}

\makeatletter
\def\cnewtheorem#1[#2]#3{\newtheorem{#1}{#3}
\expandafter\let\csname c@#1\endcsname\c@Thm}
\makeatother

\AtBeginDocument{\let\bar\wbar\let\tilde\wtilde\let\hat\what}

\newtheorem{Thm}{Theorem}
\cnewtheorem{MThm}[Thm]{Main Theorem}
\cnewtheorem{Coro}[Thm]{Corollary}
\cnewtheorem{Lem}[Thm]{Lemma}
\newtheorem{Qu}{Question}

\theoremstyle{definition}
\cnewtheorem{Def}[Thm]{Definition}

\makeautorefname{Thm}{Theorem}
\makeautorefname{MThm}{Main Theorem}
\makeautorefname{Coro}{Corollary}
\makeautorefname{Qu}{Question}
\makeautorefname{Lem}{Lemma}
\makeautorefname{Def}{Definition}

\begin{document}

\begin{abstract} 
We define integral measures of complexity for Heegaard splittings based on the
graph dual to the curve complex and on the pants complex defined by Hatcher 
and Thurston.  As the Heegaard splitting is stabilized, the sequence of 
complexities turns out to converge to a non-trivial limit depending only on the
manifold.  We then use a similar method to compare different manifolds, 
defining a distance which converges under stabilization to an integer related 
to Dehn surgeries between the two manifolds.
\end{abstract}

\maketitle

\section{Introduction}

Hempel introduced the curve complex to the study of Heegaard splittings by 
defining a distance which generalizes the definitions of reducible, weakly 
reducible and the disjoint curve property.  This has proved very useful in 
studying irreducible splittings of manifolds.  Hempel's definition of distance 
is, in some sense, very local.  If a Heegaard splitting is stabilized only 
once, the distance drops to zero, regardless of the original splitting.  

This definition of distance can be directly modified to use the pants complex 
defined by Hatcher and Thurston, or the closely related dual graph to the curve
complex.  (All three spaces will be defined in detail in the next section.) The
two types of distance that come from these metric spaces prove to give a more 
global measure of complexity.  In most cases, when a splitting is stabilized, 
the distance will increase by one (rather than dropping to zero.)  Thus if the 
genus is subtracted from the distance, the resulting integer, called the 
complexity, should tend to stay constant under stabilization.

We will show that for any Heegaard splitting of a given manifold, under
an infinite sequence of stabilizations the complexities of the stabilizations
form a convergent sequence whose limit is non-trivial and which depends only 
on the manifold.

\fullref{defsect} contains the definitions of the three spaces mentioned
above.  In \fullref{maxdisksect}, we consider maximal collections of
disks in compression bodies.  This will allow us to apply the later results
to manifolds with boundary as well as closed manifolds.  We define the
dual distance and pants distance in \fullref{distsect} and give lower
bounds for this distance depending on the ambient 3--manifold.

In \fullref{stabsect}, we define the complexity of a Heegaard splitting
and show that the sequence of complexities converges.  The limit is called the
Heegaard complexity of the manifold.  The basic properties of the Heegaard
complexity are examined in \fullref{propsect}.  

Finally, in \fullref{compmfldsect}, we show that a similar technique can 
be used to define a measure of the distance between different manifolds.  
This distance turns out to be equal to the minimal number of components needed 
for a link such that Dehn surgery on the link in one manifold produces the
second manifold.  \fullref{questsect} is a list of questions that
arise and speculations about applications.

I would like to thank my advisor, Abby Thompson, for her guidance and
support.  This research was supported by NSF VIGRE grant 0135345.

\section{Definitions}
\label{defsect}

Let $\Sigma$ be a compact, connected, closed, orientable surface.  Throughout
the paper, we will assume the genus of $\Sigma$ is greater than 1.
We will give the definitions of three metric spaces based on $\Sigma$.

The first, the \emph{curve complex}, $C(\Sigma)$, is the cell complex
defined as follows:  The vertices of $C(\Sigma)$ are
isotopy classes of non-trivial, simple closed curves in $\Sigma$.
An edge will connect two vertices if and only if there are
representatives of the two isotopy classes which are disjoint.  

The graph can be filled in with cells of higher dimension. A
collection of vertices $\{u_0, \dots, u_n\}$ bounds an $n$--simplex in
$C(\Sigma)$ if and only if some collection of loops defined by the 
vertices are pairwise disjoint. Let $g$ be the genus of $\Sigma$. Then 
$C(\Sigma)$ has dimension $3g - 4$, a maximal simplex corresponding to a
pair-of-pants decomposition of $\Sigma$, with $3g - 3$ loops.

The second space, the \emph{dual curve complex}, $C^*(\Sigma)$,
will be defined as follows: each vertex $v \in C^*(\Sigma)$
corresponds to a maximal dimensional simplex $\sigma_v$ in
$C(\Sigma)$, ie a pair-of-pants decomposition of $\Sigma$.
We will use the convention that $u$ is a vertex of $C(\Sigma)$
and $v$ is a vertex of $C^*(\Sigma)$.

Two vertices, $v$, $v'$ are connected by an edge if
and only if the simplices $\sigma_v$ and $\sigma_{v'}$ in
$C(\Sigma)$ share a co-dimension one face.
An edge thus corresponds to a move of the following type:
Start with a pants decomposition of $\Sigma$; remove one loop and
replace it with a loop which is disjoint from all the other loops
and which creates a new pants decomposition.

The third space, the \emph{pants complex}, $C^P(\Sigma)$, was
defined by Hatcher and Thurston~\cite{HT:rep} as a discrete quotient
of the space of Morse functions on $\Sigma$.  Every Morse function
on $\Sigma$ suggests a decomposition of $\Sigma$ into pairs of
pants, forming the vertices of $C^P(\Sigma)$.  Near-Morse
functions, in which there are two critical points at the same
level, form co-dimension one sets between the connected components
of Morse functions.  These suggest edges between certain vertices.

Careful analysis shows that the edges correspond to moves of the following 
type: Given a pair-of-pants decomposition of $\Sigma$, remove a loop
from the pants decomposition and replace it with a loop which is
disjoint from the other loops and intersects the original loop
minimally.  A minimal intersection between two loops is defined
as follows:
\begin{figure}[htb]
  \begin{center}
  \includegraphics[width=3.5in]{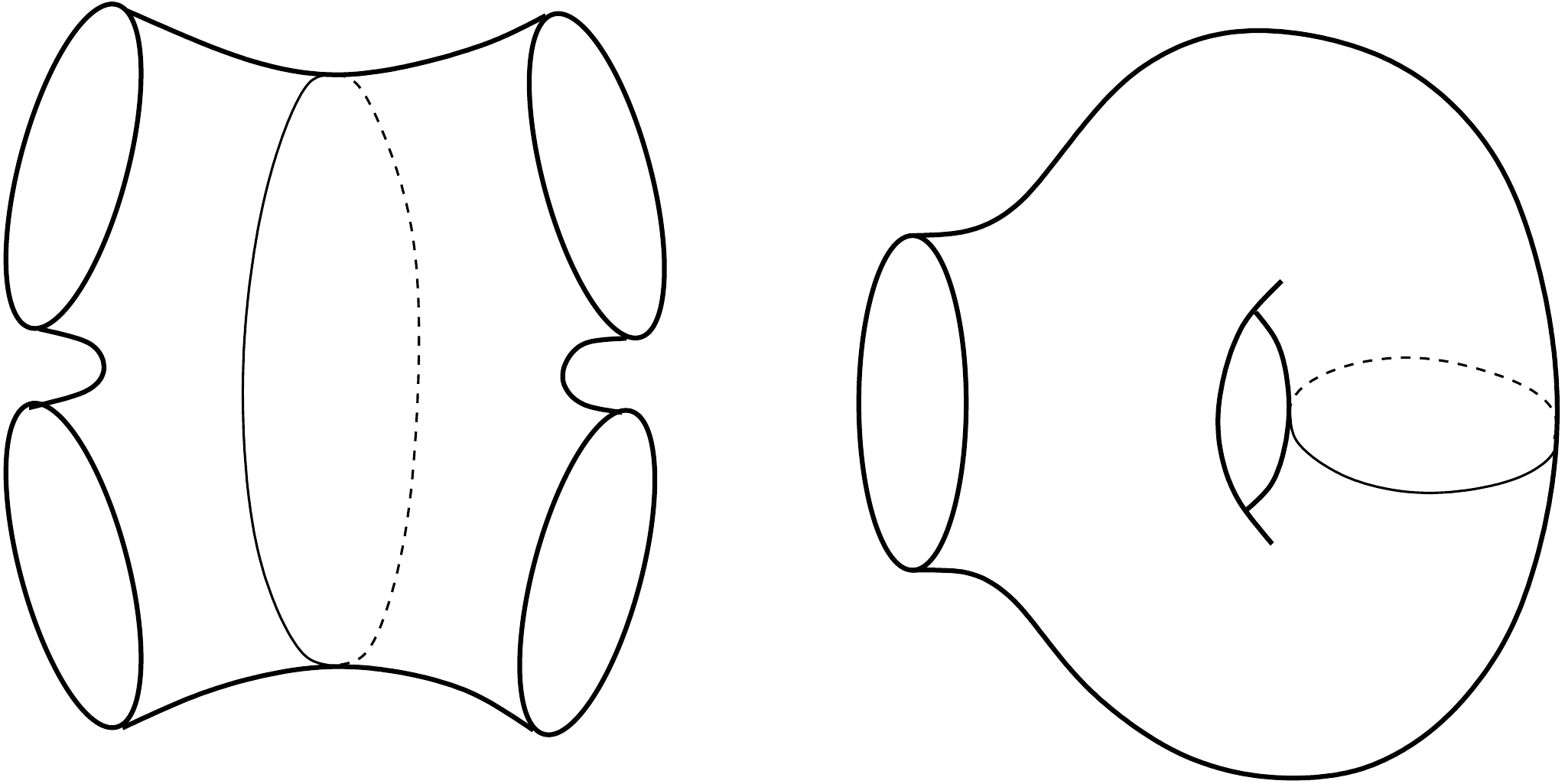}
  \caption{The two types of loops in a pair-of-pants decomposition}
  \label{loopcomps}
  \end{center}
\end{figure}
Consider a pair-of-pants decomposition $\mathbf{L} = \{l_0,\dots,l_m\}$.
The complement of $\mathbf{L}$ in $\Sigma$ is a collection of pairs of
pants.  The complement of $\{l_0,\dots, \hat l_i,\dots,l_m\}$ contains
a number of pairs of pants and one component, containing $l_i$, which is
not a pair of pants.  Call this component $C$. The surface $C$ is the result 
of either gluing two pairs together along the loop $l_i$ or of gluing 
together two cuffs of the same pair.  In the first case, $C$ is a
sphere with four punctures.  In the second 
case, $C$ is a torus with a single puncture. (See \fullref{loopcomps}.)

A loop $l'_i$, disjoint from and not parallel to 
$l_0,\dots,\hat l_i,\dots,l_m$ must lie in $C$.  If $C$ is a four-punctured
sphere then $l'_i$ intersects $l_i$ minimally when there are two points
of intersection. (The loop is separating so there must be an even 
number of intersections.) If $C$ is a once-punctured torus, then one point 
of intersection is minimal.

Because the space of Morse functions on $\Sigma$ is connected,
Hatcher and Thurston were able to
show that the pants complex is connected.  There is a canonical
map from the vertices of $C^P(\Sigma)$ to the vertices of
$C^*(\Sigma)$.  For each edge of $C^P(\Sigma)$, there is a
corresponding edge of $C^*(\Sigma)$ (but not vice versa) so there
is a canonical embedding of $C^P(\Sigma)$ into $C^*(\Sigma)$ which
is onto the vertices of $C^*(\Sigma)$.  Since $C^P(\Sigma)$ is
connected, it follows that $C^*(\Sigma)$ is connected.

Brock~\cite{Brock:pants} has shown that the pants complex, $C^P(\Sigma)$, is
quasi-isometric to Teichmuller space with the Weil--Peterson metric
and that distances in the pants complex are related to the volumes
of convex cores of hyperbolic manifolds of the form $\Sigma \times
\mathbb{R}$.  We will study the pants complex and its relative,
the dual curve complex, in relation to Heegaard splittings, in
analogue with Hempel's work with the curve complex~\cite{Hempel:complex}.

Let $H$ be a handlebody and let $\phi \co \Sigma \rightarrow
\partial H$ be a homeomorphism.  For a vertex $u \in C(\Sigma)$,
write $u \in H$ if for some loop $l$ in the isotopy class
corresponding to $u$, $\phi(l)$ bounds a disk in $H$.  (Note: If
this is true for one loop in the isotopy class, then it is true
for all the loops in the isotopy class.)

An edge path between two vertices $u$, $u'$ in the curve complex corresponds
to a sequence of loops $u \ni l_0, \dots, l_n \in u'$ in $\Sigma$ such that
consecutive loops are disjoint.  The length of the path is $n$ and the
distance $d(u,u')$ is the length of the shortest possible path.  This is
often called the \emph{geodesic metric} on $C(\Sigma)$.

Let $(\Sigma, H_1, H_2)$ be a Heegaard splitting of a manifold $M$.  
Consider the inclusion maps $\Sigma \rightarrow H_i$.  Each map suggests a set 
of vertices in $C(\Sigma)$ which are in $H_i$.  The standard distance of 
$\Sigma$, as in Hempel~\cite{Hempel:complex}, is 
$d(\Sigma) = \min\{d(u,u') | u \in H_1, u' \in H_2\}$.
This distance measures the irreducibility of $\Sigma$, in the sense that
if $d(\Sigma) = 0$ then $\Sigma$ is reducible, if $d(\Sigma) = 1$
then $\Sigma$ is weakly reducible and if $d(\Sigma) = 2$ then $\Sigma$
has the disjoint curve property.  

There are analogous definitions for the distance of a Heegaard splitting
based on the dual curve complex and the pants complex.  We will define
these right after a short aside about collections of disks in compression
bodies.

\section{Maximal collections of disks}
\label{maxdisksect}

A 1--handle is ball parameterized as $D \times [0,1]$, where $D$ is a disk. 
Let $F$ be a surface (not necessarily connected) with no sphere components.  A 
compression body is either a handlebody or a connected manifold constructed by 
gluing 1--handles to the boundary component $F \times \{1\}$ of 
$F \times [0,1]$ along the disks $D \times \{0\}$ and $D \times \{1\}$.  For a 
compression body $H$, let $\partial_- H$ = $F \times \{0\}$ and let
$\partial_+ H = \partial H \setminus \partial_- H$.
If $H$ is a handlebody then $\partial_+ H = \partial H$ and 
$\partial_- H = \emptyset$.

Let $H$ be a compression body and let $\mathbf{D} = \{D_0,\dots,D_n\}$ be a 
collection of pairwise disjoint, pairwise non-parallel, properly embedded disks
in $H$.

\begin{Def}
The collection $\mathbf{D}$ is \emph{maximal} if any properly embedded disk
in $H$ that is disjoint from $\mathbf{D}$ is parallel to one of the disks
in $\mathbf{D}$.
\end{Def}

If $H$ is a handlebody, a maximal collection of disks gives a 
pair-of-pants decomposition for $\partial H$.  Otherwise, when 
$\partial_- H \neq \emptyset$, there is no collection of disks which 
cuts $\partial H$ into pairs of pants, so we have to be more careful.  A 
maximal collection of disks will cut $H$ into a collection of balls and 
pieces homeomorphic to $\partial_- H \times I$.  However, not every 
collection of disks with this property will be maximal. 

\begin{Lem}
\label{cutlem}
Let $\mathbf{D}$ be a collection of disjoint, properly embedded disks in $H$
(not necessarily maximal).
The closure of the complement of $\mathbf{D}$ is a collection of balls and 
compression bodies.
\end{Lem}

The proof of this Lemma is left to the reader.  Let $H'$ be a component
of $H \setminus \mathbf{D}$ and let $\wwbar{H}'$ be its closure.  Then
$\wwbar{H}' \setminus H'$ is a collection of disks along which $H$ was cut.
We will call a component of $\wwbar{H}' \setminus H'$ a \emph{scar} on $H'$.

\begin{Lem}
Let $H$ be a compression body and not a ball, a solid torus or a 
$(\textit{surface}) \times I$ and let 
$\mathbf{D}$ be a  maximal collection of disjoint, non-parallel essential
disks. Then each component of $H \setminus \mathbf{D}$
is either a ball with three scars or a piece of the form $(\textit{surface})
 \times I$ with exactly one scar.
\end{Lem}

\begin{proof}
Let $H'$ be a component of $H \setminus \mathbf{D}$.  From \fullref{cutlem},
we know that $H'$ is either a ball or a compression body.

First assume $H'$ is a ball.  If there is one scar on $H'$ then this scar
corresponds to a boundary parallel disk in $\mathbf{D}$, but we assumed all 
the disks in $\mathbf{D}$ were essential, so this is impossible.  If there
are two scars on $H'$ then there are two disks in $\mathbf{D}$ which are
parallel in $H$, and again we assumed this is not the case.

If there are more than three scars on $H'$ then $\partial \wwbar{H}' \cap H'$
is a sphere with at least four punctures.  Let $l$ be an essential loop
in this surface which is not boundary-parallel.  Then $l$ bounds a disk in
$H'$, and therefore in $H$.  This disk is properly embedded, disjoint from
the rest of the disks, and not parallel to any of them.  Thus the maximality
assumption implies that every ball component $H'$ must have exactly three
scars.

Now assume $\wwbar{H}'$ is a compression body.  By definition, $H'$ is the
result of attaching 1--handles to components of the form $(\textit{surface})
 \times I$.  A 1--handle would define an essential disk which is not parallel 
to any of the disks in $\mathbf{D}$, so maximality implies there cannot be any
1--handles and $\wwbar{H}'$ must be of the form $F \times I$ where $F$ is a
surface and not a sphere.

Because $H$ is connected and $H'$ is not all of $H$, there must be at least
one scar on $H'$.  If there are two or more scars on $H'$, let $l_1$, $l_2$
be the boundaries of two scars in $\partial \wwbar{H}$.  Let $\alpha$ be
an arc from $l_1$ to $l_2$, disjoint from the rest of the scars, and let $N$ 
be a regular neighborhood of $l_1 \cup l_2 \cup \alpha$ in $\partial \wwbar{H}'$. 
The disk $\partial N \setminus \partial \wwbar{H}$ will be disjoint from the scars
in $\partial \wwbar{H}'$ and boundary--parallel in $H'$.  Because of the scars 
bounded by $l_1$ and $l_2$, this disk is non-trivial in $H$ and not parallel to
any of the disks in $\mathbf{D}$.  Thus maximality of $\mathbf{D}$ implies that
$H'$ must have exactly one scar.
\end{proof}

\begin{Lem}
\label{nonseplem1}
Let $\Sigma$ be a positive-genus surface and let $\mathbf{L}$ be a 
pair-of-pants decomposition for $\Sigma$.  Then some loop 
$l \in \mathbf{L}$ is non-separating.
\end{Lem}

\begin{proof}
Cutting a positive-genus surface along a separating loop produces two
surfaces (with boundary), each with strictly positive genus.  By induction, if 
we cut $\Sigma$ along all the separating loops in $\mathbf{L}$, the result
will be a number of positive-genus surfaces.  Since $\mathbf{L}$ is a
pair-of-pants decomposition, cutting along all the loops should produce
a collection of pairs of pants (genus-zero surfaces) so at least one of
the loops in $\mathbf{L}$ must be non-separating.
\end{proof}

\begin{Lem}
\label{nonseplem2}
Let $\Sigma$ be a closed surface of genus $g$.  Let $\mathbf{L}$ be a 
pair-of-pants decomposition of $\Sigma$ and
let $l_1,\dots,l_n \in \mathbf{L}$ be distinct loops such that their union
does not separate $\Sigma$.  Then there are loops
$l_{n+1},\dots,l_g \in \mathbf{L}$ such that $l_1,\dots,l_g$ are distinct
loops whose union does not separate $\Sigma$.
\end{Lem}

\begin{proof}
Given loops $l_1,\dots,l_n$, we will show that if $n < g$ then there
is a distinct loop $l_{n+1} \in \mathbf{L}$ such that the union of 
$l_1,\dots,l_{n+1}$ does not separate $\Sigma$.  By taking the collection to be
maximal, this implies the desired result.

Let $\Sigma'$ be the result of cutting $\Sigma$ along the loops 
$l_1,\dots,l_n$ and gluing a disk into each boundary component of the 
resulting surface.  Because the union of the loops $l_1,\dots,l_n$ does
not separate $\Sigma$, the surface $\Sigma'$ is connected.  We also assumed
$n < g$, so $\Sigma'$ has strictly positive genus.

The remaining loops in $\mathbf{L}$ contain a pair-of-pants decomposition
for $\Sigma'$. (Some of the loops remaining in $\mathbf{L}$ will be trivial
or parallel in $\Sigma'$ and need to be thrown out.)  By 
\fullref{nonseplem1}, there is a loop $l_{n+1} \in \mathbf{L}$ in the
induced pair-of-pants decomposition which does not separate $\Sigma'$.
Then $l_1,\dots,l_{n+1}$ does not separate $\Sigma$. By induction, the proof is
complete.
\end{proof}

Given a compression body $H$, let $g$ be the genus of $\partial_+ H$ and
let $b$ be the sum of the genera of the components of $\partial_- H$.

\begin{Lem}
\label{nonseplem3}
Let $\mathbf{D}$ be a maximal collection of disks for $H$ and let 
$\mathbf{L}$ be a pair-of-pants decomposition for $\partial_+ H$ such
that for each disk $D \in \mathbf{D}$, $\partial D \in \mathbf{L}$.
If $l_1,\dots,l_g \in \mathbf{L}$ is a collection of distinct loops whose 
union does not separate $\Sigma$, then at least $g - b$ of these loops
bound disks in $\mathbf{D}$.
\end{Lem}

\begin{proof}
Cut $H$ along the disks in $\mathbf{D}$.  Let $H'$ be a resulting component
of the form $F \times I$ where $F$ is a surface of genus $g'$.  Any 
collection of more than $g'$ loops in $F$ will separate $F$, so if more than
$g'$ of the loops $l_1,\dots,l_g$ are in $H'$, then these loops will separate
$\partial_+ H'$.  Because there is exactly one scar on $H'$, these disks 
would also separate $H$.  Since we assumed that the union of $l_1,\dots,l_g$
does not separate $H$, at most $g'$ of these loops can be in $H'$.

The same is true for any non-ball component of $H \setminus \mathbf{D}$,
so the number of loops that are not boundaries of disks is at most $b$.
\end{proof}

\section{Distance and Heegaard splittings}
\label{distsect}

We will now define a notion of distance for Heegaard splittings based on 
the dual curve complex.
Let $v$ be a vertex of $C^*(\Sigma)$.  We will say $v$ \emph{defines} a 
compression body $H$ if there is a maximal collection of disks for $H$ such 
that the boundary of each disk defines a vertex $u$ of $\sigma_v$.

\begin{Lem}
Assume $v$ defines a compression body $H$, with 
$\phi \co \Sigma \rightarrow \partial_+ H$ and a second compression body $H'$
with $\phi'\co \Sigma \rightarrow \partial_+ H'$.  If the same vertices of $v$ 
bound disks in $H$ as in $H'$ then there is a homeomorphism 
$\psi \co H \rightarrow H'$ such that $\psi \circ \phi = \phi'$.
\end{Lem}

A proof of the Lemma is left to the reader. The converse is not true.  In 
general, given two pants decompositions of $\Sigma$, there may not be an 
automorphism of $\Sigma$ taking one to the other.  For example, if one of the 
pants decompositions consists entirely of non-separating loops (such a 
decomposition exists) then there will be no homeomorphism taking it to a pants 
decomposition containing separating loops.

The \emph{dual distance} $D(v, v')$ between two vertices in
$C^*(\Sigma)$ is the length of the shortest path in $C^*(\Sigma)$ between
them.  For a Heegaard splitting, $(\Sigma, H_1, H_2)$, the 
\emph{dual distance} of $\Sigma$ is
$D(\Sigma) = \min\{D(v,v') | v$ defines $H_1, v'$ defines $H_2\}$.
Note that $D(\Sigma) \geq d(\Sigma)$.  Hempel has shown that there are
genus two Heegaard splittings such that $d(\Sigma)$ is arbitrarily large.  
Thus there are Heegaard splittings with $D(\Sigma)$ arbitrarily large.

Let $D^P(v,v')$ be the distance between vertices $v$ and $v'$ in the
pants complex.  Because of the one-to-one map between the vertices of
$C^*(\Sigma)$ and the vertices of $C^P(\Sigma)$, we can think of $v$ and
$v'$ as being in either graph.  An edge path in $C^P$ maps to an edge path of
the same length in $C^*$ so $D(v,v') \leq D^P(v,v')$.

Let $(\Sigma, H_1, H_2)$ be a genus $g$ Heegaard splitting of a 3--manifold
$M$.  From now on we will assume that $H_2$ is a handlebody (ie 
$\partial_- H_2 = \emptyset$) but we will allow $H_1$ to be a compression
body.  Thus $\partial_- H_1 = \partial M$ and $H_1$ will be a handlebody
if and only if $M$ is closed.  Such a Heegaard splitting always exists.

Let $b$ be the sum of the genera of the boundary components of $M$ and let $n$ 
be the maximal number of disjoint, embedded 2--spheres $S_1,\ldots,S_n$ such 
that $M \setminus \bigl(\bigcup S_i\bigr)$ is connected.  (Equivalently, $n$ is the 
number of $S^1 \times S^2$ components of the prime decomposition of $M$, so $n$
is well defined and finite.)

\begin{Lem}
\label{negblem1}
$D(\Sigma) \geq g - b - n$
\end{Lem}

\begin{proof}
Let $v$ define $H_1$ and $v'$ define $H_2$ and assume for contradiction
$D(v,v') = D(\Sigma) < g - b - n$.  Let $\mathbf{L}$ be a collection of 
pairwise-disjoint loops in $\Sigma_g$ corresponding to the vertices of 
$\sigma_v$.  Recall that $H_1$ may be a compression body (if $b > 0$) or a 
handlebody (if $b = 0$).

By \fullref{nonseplem2}, there are loops $l_1,\dots,l_g \in \mathbf{L}$
whose union does not separate $\Sigma$.  Each step in the path from $v$ to
$v'$ changes one vertex.  Since $D(v,v') \leq g - (b + n + 1)$ and there are 
$g$ loops in the collection, at least $b + n + 1$ of the loops $l_1,\ldots,l_g$
correspond to vertices of $\sigma_{v'}$.  Since $H_2$ is a handlebody,
these $b+n+1$ loops bound disks in $H_2$.  
By \fullref{nonseplem3}, at most $b$ of the loops do not bound disks
in $H_1$ so at least $n+1$ of the loops $l_1,\dots,l_2$ bound disks in both
$H_1$ and $H_2$.  

Let $l_1,\dots,l_{m}$ be the loops that bound disks in
both $H_1$ and $H_2$.  Then each $l_i$ defines an embedded sphere $S_i$.
These spheres are disjoint and because the loops $l_1,\dots,l_m$ do not
separate $\Sigma$, the spheres $S_1,\dots,S_m$ do not separate $M$ (hence
they are non-parallel).
We assumed that $M$ contains at most $n$ such spheres, so we must have 
$m \leq n$.  However, we showed that there are at least $n+1$ loops bounding
disks on both sides.  This contradiction completes the proof.
\end{proof}

\begin{Coro}
If $M$ is closed and irreducible, then $D(\Sigma) \geq g$.
\end{Coro}

\begin{Lem}
\label{negblem2}
$D^P(\Sigma) \geq g - b - n$
\end{Lem}

The proof of \fullref{negblem2} is identical to the proof of 
\fullref{negblem1}, after replacing each $D$ with $D^P$.

\begin{Lem}
\label{negblem3}
If $D(\Sigma) = g - b - n$ then $M = S^3$ or $M$ is a connect sum of
lens spaces, handlebodies and copies of $S^1 \times S^2$.
\end{Lem}

\begin{proof}
Let $v, v'$ be vertices of $C^*(\Sigma)$ that define $H_1$ and $H_2$,
respectively, so that $D(v,v') = D(\Sigma) = g - b - n$.
As in the previous proof, let $\mathbf{L}$ be the collection of loops 
corresponding to the vertices of $\sigma_v$ and let 
$l_1,\dots,l_g \in \mathbf{L}$ be a collection of loops whose union is
non-separating. 

We saw that if more than $b + n$ of these loops also 
correspond to vertices of $\sigma_{v'}$ then we have a contradiction.  
Because $D(v,v') = g - b - n$, and there are $g$ loops, at least $b + n$
of them must be common to $v$ and $v'$ so we know that exactly $b + n$
of the loops $l_1,\dots, l_g$ are common to both $\sigma_v$ and $\sigma_{v'}$
and each of the remaining loops is moved exactly once.  The remaining loops
of $\mathbf{L}$ are not moved.

Of the $b + n$ loops that are common to $\sigma_v$ and $\sigma_{v'}$,
at most $n$ bound disks in $H_1$ (because there are at most $n$ 
non-separating spheres in $M$) so at least $b$ of the loops do not bound
disks in $H_1$.  Thus each of the loops in $l_1,\dots,l_g$ that is moved
must bound a disk in $H_1$.

Assume $l_1$ is the first loop that is moved, and replaced by a loop $l_1'$.
Since $H_2$ is a handlebody and the loop $l_1'$ is not moved later in
the path in $C^*(\Sigma)$, $l_1'$ must bound a disk in $H_2$.  Recall that
there are two types of loops in the pants decomposition, defined by whether
removing the loop from the collection produces a four-punctured sphere or
a once-punctured torus.

If $l_1$ is of the first type, then the four-punctured sphere is part of
the boundary of a ball in $H_1$.  Since $l_1'$ sits in the four-punctured
sphere, it bounds a disk in the ball in $H_1$, and therefore bounds a disk
in $H_1$.  Thus we could have started the path with $l_1'$ in the pants
decomposition instead of $l_1$.  Because the path in $C^*(\Sigma)$ was
assumed to be minimal, $l_1$ must sit in a punctured torus.

Let $l''_1$ be the loop defining
the boundary of the puncture.  Notice that $l''_1$ is separating, so it cannot
be one of the non-separating loops $l_1,\dots,l_g$ and it cannot be moved
later on.  In particular, this
implies that $l_1$ cannot be adjacent to any of the loops that are moved
later in the path.

The loop $l''_1$ bounds a disk in $H_1$ and a disk in $H_2$, so $l''_1$ 
defines a sphere in $M$ which cuts off a genus-one piece of the Heegaard 
splitting.  Thus the sphere defines either a stabilization or a connect sum 
with a lens space.

Assume $l_2$ is the next loop which is moved. We saw that $l_2$ cannot be 
adjacent to $l_1$, so $l_2$ sits in a punctured torus, with a punctured 
bounded by $l''_2$.  Again, $l''_2$
defines a sphere which separates a genus-one piece of the Heegaard splitting,
so $l''_2$ defines a stabilization or a connect sum with a lens space.
Continuing in this fashion for each loop that is moved, we get a collection
of $D(\Sigma)$ stabilizations and lens space summands.

Let $M'$ be the result of cutting $M$ along these spheres and gluing balls
into the resulting boundary components.  Let $(\Sigma', H_1', H_2')$ be the 
Heegaard splitting resulting from gluing disks into the image of $\Sigma$ 
in $M'$.  Since we removed from $M$ all the loops that were not in both the 
pants decomposition of $H_1$ and that of $H_2$, it follows that 
$D(\Sigma') = 0$.  

Let $F$ be a component of $\partial_- H_1$.  Because the pants decomposition
for $\Sigma$ was maximal, there is a loop $l''$ which bounds a disk in
$H_1$ and separates $H_1$ into a compression body and an $F \times I$
component.  Since $D(\Sigma') = 0$, the loop $l''$ bounds a disk in $H_2$ and
defines a separating sphere in $M'$.  The disks cut off a handlebody from
$H_2$ and a $F \times I$ component from $H_1$, so the sphere cuts off
a handlebody summand from $M'$.

For each boundary component of $M'$, (each corresponding to a component of
$\partial_- H_1$) there is a corresponding loop defining a sphere which cuts 
off a handlebody summand from $M$.  Let $M''$ be the result of cutting $M'$ 
along these spheres and gluing balls into the resulting spheres.  Let
$(\Sigma'', H_1'', H_2'')$ be the Heegaard splitting resulting from repairing
the image of $\Sigma'$ in $M''$.

If $H_1''$ is a ball then $M''$ is $S^3$ and we are done.
Otherwise, $M''$ is closed and $D(\Sigma'') = 0$, so $M''$ is the
result of gluing together two handlebodies by the identity map on their
boundaries.  This construction always yields a connect sum of copies of
$S^1 \times S^2$.  This completes the proof.
\end{proof}

\begin{Lem}
\label{negblem4}
If $D^P(\Sigma) = g - b - n$ then $M = S^3$ or $M$ is a connect sum of
handlebodies and copies of $S^1 \times S^2$.
\end{Lem}

Again, the proof of this lemma is almost identical to the analogous proof for
$D(\Sigma)$.  Note, however, that the set of manifolds in the second lemma
is more restricted.

\section{Stabilization}
\label{stabsect}

Let $(\Sigma, H_1, H_2)$ be a Heegaard splitting of $M$ with $b$ 
and $n$ defined as in \fullref{distsect}.  For each $g$ greater
then or equal to the genus of $\Sigma$, let $\Sigma_g$ be a stabilization of
$\Sigma$ such that $\Sigma_g$ has genus $g$.  In other words, $\Sigma_g$ is
the result of attaching zero or more trivial handles to $\Sigma$ so that the
resulting surface is a Heegaard surface.

Define $A_g(\Sigma) = D(\Sigma_g) + b - g$ and 
$A^P_g(\Sigma) = D^P(\Sigma_g) + b - g$.  We will consider the limiting
behavior of these two values.
From Lemmas \ref{negblem1} through \ref{negblem4}, the following two
Lemmas follow immediately:

\begin{Lem}
$A_g(\Sigma) \geq -n$ and if $A_g(\Sigma) = -n$ then $M = S^3$ or $M$ is 
a connect sum of handlebodies, copies of $S^1 \times S^2$ and lens spaces.
\end{Lem}

\begin{Lem}
$A^P_g(\Sigma) \geq -n$ and if $A_g(\Sigma) = -n$ then $M = S^3$ or $M$ is 
a connect sum of handlebodies and copies of $S^1 \times S^2$.
\end{Lem}

This gives us a lower bound on the sequences $A_g(\Sigma)$ and 
$A^P_g(\Sigma)$ as $g \rightarrow \infty$.  We will next show that the
sequences are also bounded above.  For the following lemmas, consider a 
fixed $g$ and a sequence of stabilizations $\Sigma_h$.

\begin{Lem}
\label{uprbndlem}
If $h > g$ then $A_h(\Sigma) \leq 2A_g(\Sigma) + g - b$
and  $A^P_h(\Sigma) \leq 2A^P_g(\Sigma) + g - b$.
\end{Lem}

\begin{proof}
We will prove that $D(\Sigma_h) \leq 2 D(\Sigma_g) + (h-g)$.  By subtracting
$h$ from both sides we get the stated result.
The pants distance, $D^P$ can be
substituted for $D$ throughout the proof.
Let $v_1,\ldots, v_n$ be a minimal path from $H_1$ to $H_2$ in $C^*(\Sigma)$
(or $C^P(\Sigma)$).
Let $l_1^1,\dots,l_m^1$ ($m = 3g -3$) be loops in $\Sigma$ corresponding to
the vertices of $\sigma_{v_1}$.

For each $v_i$, let $l^1_i,\dots,l^m_i$ be loops corresponding to the vertices
of $\sigma_{v_i}$ and assume the loops are labeled so that if $l^j_{i-1}$ is 
a vertex of
both $\sigma_{v_{i-1}}$ and $\sigma_{v_i}$ then $l^j_i = l^j_{i-1}$.  In other
words, if the move from $v_{i-1}$ to $v_i$ replaces a loop $l_{i-1}^j$ with
a new loop, this new loop is labeled $l_i^j$.  If a loop $l_{i-1}^j$ is
not replaced, the same loop appears in $v_i$ as $l_i^j$.

\begin{figure}[htb]
  \begin{center}
  \includegraphics[width=3.5in]{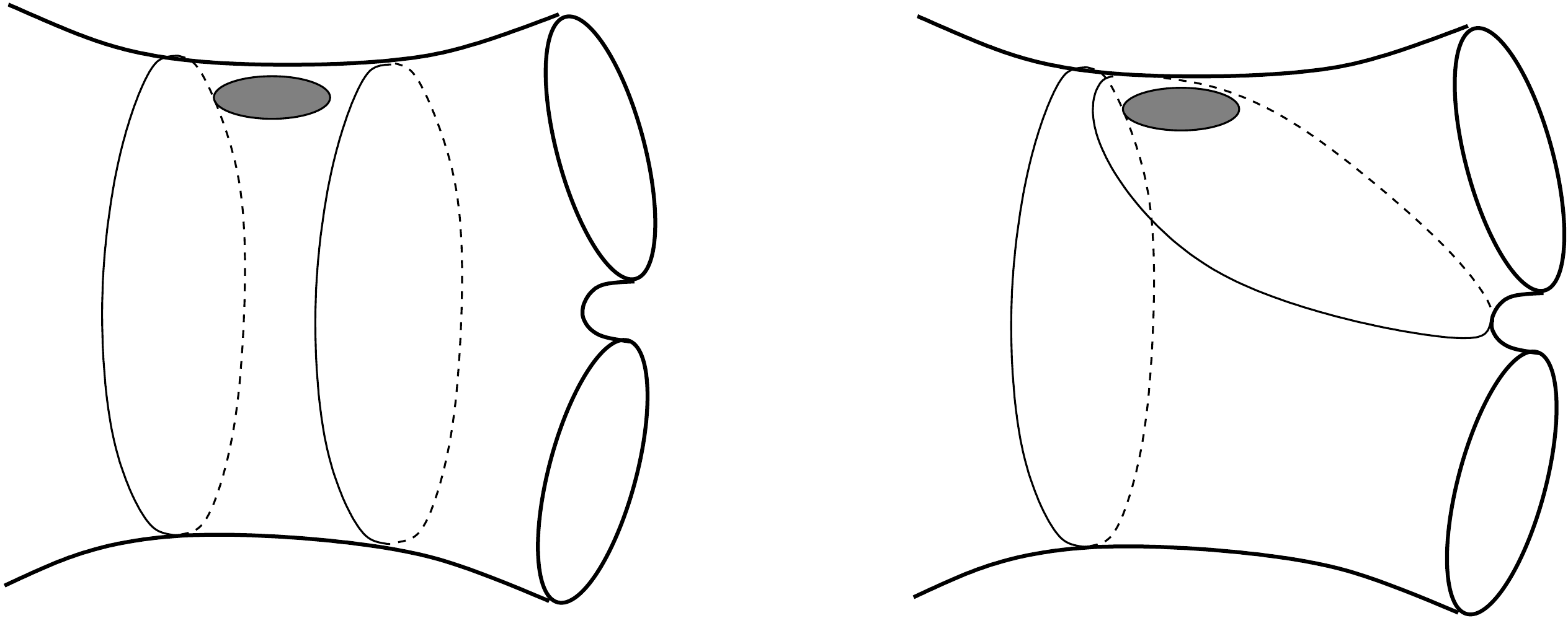}
  \caption{One extra move is required to get the extra loop out of the
     way.}
  \label{extramove}
  \end{center}
\end{figure}
We will define a stabilization of $\Sigma_g$ and a pants decomposition
$k^1_1,\dots,k^{m'}_1$ ($m'=3h-3$) as follows:  Let $k^j_1 = l^j_1$
for $1 \leq j \leq m$.  Let $k^{m+1}_1$
be a loop parallel to $k^1_1$ and let $k^{m+2}_1$ be a trivial loop in the
resulting annulus.  See \fullref{extramove}.

\begin{figure}[htb]
  \begin{center}
  \includegraphics[width=3.5in]{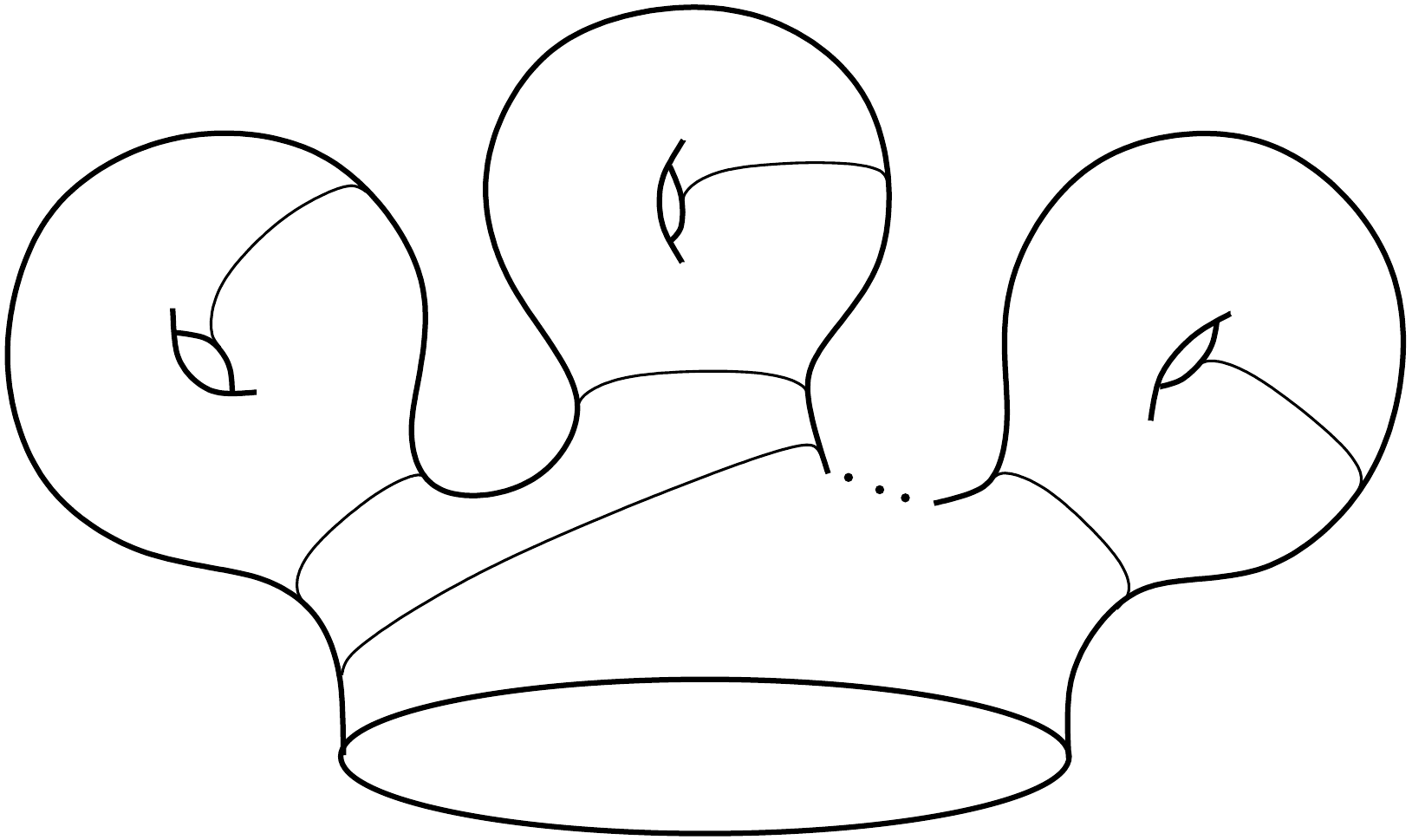}
  \caption{A stabilization glued into $\Sigma$}
  \label{stabloops}
  \end{center}
\end{figure}
Puncture $\Sigma$ in the disk defined by $k^{m+2}_1$
and construct a new surface $\Sigma'$ by attaching a punctured genus--$(h-g)$
surface to $\Sigma$ (so that the resulting surface is a stabilization).  Let 
$k^{m+3}_1,\dots,k^{m'}_1$ be loops on $\Sigma'$
as shown in \fullref{stabloops} and let $(\Sigma', H_1', H_2')$ be the 
resulting Heegaard splitting.  Note that exactly $(h - g)$ of the new
loops are non-separating and the rest of the loops bound disks in both
$H_1'$ and $H_2'$.
The vertex $v'_1 \in C^*(\Sigma')$ given by 
the loops $k^1_1,\dots,k^{m'}_1$ defines the compression body $H_1'$.

We will create a sequence of vertices $v'_1,\dots,v'_{n''}$ in $C^*(\Sigma')$
such that $n'' \leq 2n + (h-g)$ and $v'_{n''}$ defines $H_2'$.
If $l_1^1 = l^1_2$ then for $j \leq m$, let $k^j_2$ be the loop in $\Sigma'$
defined by $l^j_2$ in $\Sigma$ and for $j > m$, let
$k^j_2 = k^j_1$.  One can check that the vertex $v'_2 \in C^*(\Sigma')$ is 
connected to $v'_1$ by an edge.

If $l^1_1 \neq l^1_2$, then this construction does not work because $l^1_2$
will intersect $k^{m+1}_1$.  We need to get
$k^{m+1}_1$ out of the way first.  Let $k^j_2 = k^j_1$ for each $j \neq m+1$.
The loop $k^{m+1}_1$ sits in a four-punctured sphere.  Two of the punctures
come from $l_1^1$ and $k^{m+2}_2$.  Let $l^j_1$ be the loop defining one
of the two remaining punctures.
Let $k^{m+1}_2$ be a loop parallel to $l^j_1$ so that $k^{m+2}_2$
is in the annulus defined by $k^{m+1}_2$ and $l^j_1$.  (See 
\fullref{extramove}.)
The vertex $v'_2$ defined by these loops is an edge away from $v'_1$.

Now we can define $k^j_3$ almost as we defined $k^j_2$ in the original case. 
Let $k^j_3 = l^j_2$ for $1 < j \leq m$ and $k^j_3 = k^j_2$ for $j > m$.  We 
cannot necessarily choose $k^1_3$ to be equal to $l^1_2$ because $l^1_2$ may 
intersect $k^{m+1}_2$.  However, the image in $\Sigma$ of $k^{m+1}_2$ is 
parallel to a loop which is disjoint from $l^1_2$.  We can let $l^1_3$ be a 
loop in $\Sigma'$ which is disjoint from $k^{m+1}_2$ and whose image in 
$\Sigma$ is isotopic to $l^1_2$.  The vertex $v'_3$ defined by these loops is 
an edge away from $v'_2$.

The vertex $v'_2$ or $v'_3$ defined in this way has the property that if we 
surger $\Sigma'$ along $k^{m+2}_2$, the loops in the component isotopic to 
$\Sigma$ define the vertex $v_2$.  Thus we can repeat the construction for 
$v_3$ through $v_n$.  At each stage, if the extra loop is parallel to a loop 
that needs to be moved, it takes one extra move to push it out of the way.  The
resulting sequence $v'_1,\dots,v'_{n'}$ will be at most twice as long as the 
original.

There are only $h - g$ moves left before $v'_i$ defines the handlebody $H_2'$.
Exactly  $h-g$ of the loops $k_{n'}^{m+3},\dots, k_{n'}^{m'}$ do not bound 
disks in $H_2$.  However, each can be replaced by a new loop, disjoint from 
the rest, 
which does bound a disk in $H_2$.  Thus in $h-g$ moves, we can complete the
sequence $v'_1,\dots v'_{n''}$ so that each loop
$k^{m+3}_{n''},\dots,k^{m+3}_{n''}$ bounds a disk in $H_2$.  By 
construction, the loop
$k^{m+2}_{n''} = k^{m+2}_1$ bounds a disk in $H_2$ and the loop
$k^{m+1}_{n''}$ is parallel to some $l^j_{n}$ in $\Sigma$ so it bounds a disk
in $H_2$.  Thus $v'_{n'}$ defines $H'_2$.
\end{proof}

By carefully choosing the sequence of loops $l_j^i$ to which $k_{m+1}^{i'}$
is parallel, one could improve the bound, but for our purposes, the existence
of a bound is all that is necessary.  The bound is true
for every $h$, but with $g$ fixed, so the sequence $A_h(\Sigma)$ is bounded
above and below as $h \rightarrow \infty$.  We will show that the sequence
actually converges.

\begin{Lem}
\label{scndbndlem}
For sufficiently large $g$, $A_h(\Sigma) \leq A_g(\Sigma)$ 
and $A^P_h(\Sigma) \leq A^P_g(\Sigma)$ whenever $h \geq g$.
\end{Lem}

\begin{proof}
Let $g'$ be the genus of $\Sigma$.
By \fullref{uprbndlem}, $A_g(\Sigma)$ is bounded by
$2 A_{g'}(\Sigma)+(g'-b)$.
Choose $g$ so that $4g - 3 > 2 A_{g'}(\Sigma) + (g' - b)$.  Then
$A_g(\Sigma) < 4g - 3$ so $D(\Sigma_g) < 3g - 3$.
This implies that for the minimal path from $H_1$ to $H_2$ in $C^*(\Sigma_g)$,
there is some loop that is not moved.

Consider the proof of \fullref{uprbndlem}.  If we had chosen $l_1$ to be
a loop that is never moved in the sequence $v_1,\dots,v_n$, then for the path
$v'_1,\dots,v'_{n'}$, we would have $n' = n + (h-g)$.  For sufficiently large 
g, we can find such a loop, so
we have $D(\Sigma_h) \leq D(\Sigma_g) + (h-g)$.  This proves the lemma.
\end{proof}

\begin{MThm}
The limits $\lim_{g \rightarrow \infty} A_g(\Sigma)$ and 
$\lim_{g \rightarrow \infty} A^P_g(\Sigma)$ exist and depend only
on $M$, not on the choice of $\Sigma$.
\end{MThm}

\begin{proof}
The sequence $A_g(\Sigma)$ is bounded below by \fullref{negblem1} and
non-increasing for sufficiently large $g$ by \fullref{scndbndlem}, so
the limit exists.

Given two Heegaard surfaces $\Sigma$ and $\Sigma'$ of $M$, there is a common
stabilization.  In other words, there is a genus $g$ and a Heegaard surface 
$\Sigma''$ such that $\Sigma''$ is isotopic to $\Sigma_g$ and $\Sigma'_g$.
Then $A_h(\Sigma) = A_h(\Sigma'') = A_h(\Sigma')$ for $h \geq g$ and the limits
are the same for $A_h(\Sigma)$ and $A_h(\Sigma')$.
\end{proof}

\begin{Def}
The \emph{Heegaard complexity} of $M$ is
$A(M) = \lim_{h \rightarrow \infty} A_h(\Sigma)$ where $\Sigma$ is
any Heegaard splitting of $M$.  The \emph{pants complexity} of $M$
is $A^P(M) = \lim_{h \rightarrow \infty} A^P_h(\Sigma)$.
\end{Def}

The sequence $A_h(\Sigma)$ consists entirely of integers so
for any $\Sigma$, $A_h(\Sigma) = A(M)$ for some $h$.
Thus if $A(M) = -n$ then $M = S^3$ or $M$ is a connect sum of compression
bodies, lens spaces and copies of $S^1 \times S^3$.

\section{Properties}
\label{propsect}

\begin{Lem}
Let $M$ and $M'$ be compact manifolds.  Then  $A(M \# M') \leq A(M) + A(M')$
and $A^P(M \# M') \leq A^P(M) + A^P(M')$
\end{Lem}

(Here, $M \# M'$ is the connect sum of $M$ and $M'$.)

\begin{proof}
Let $\Sigma$ be a Heegaard splitting of $M$.  Let $\Sigma_g$ be a
stabilization such that $D(\Sigma_g) + b - g = A(M)$ and for a minimal 
path from $H_1$ to $H_2$ in
$C^*(\Sigma_g)$, there is a loop $l_1$ which is a vertex of each
$\sigma_{v_i}$.  Let $\Sigma'_{g'}$
be a similar Heegaard splitting for $M'$ and let $l'_1$ be the loop that is
fixed.

We will construct a Heegaard splitting for $M \# M'$ as follows:
Let $l_2$ be a loop in $\Sigma_g$ parallel to $l_1$, let $B$ be a ball in
$M$ such that $M \cap \Sigma_g$ is a disk in the annulus defined by
$l_1$ and $l_2$ and define $l_3 = \partial(B \cap \Sigma_g)$.  Define
$l'_2$, $l'_3$ and $B'$ similarly in $M'$.

Take the connect sum of $M$ and $M'$ by removing $B$ and $B'$ from $M$ and
$M'$ respectively, then gluing together the resulting boundaries.  Choose
a gluing map that sends $l_3$ to $l'_3$.  Let $\Sigma''$ be the resulting
Heegaard splitting.  The paths in $C^*(\Sigma_g)$ and $C^*(\Sigma'_g)$
define a path in $C^*(\Sigma'')$, implying that.
$D(\Sigma'') \leq D(\Sigma_g) + D(\Sigma_{g'})$ and the genus of $\Sigma''$
is $g + g'$.

The sequence of stabilizations $\Sigma''_h$ can be constructed by the above
gluing operation, by taking stabilizations of $\Sigma_g$.  Thus
$D(\Sigma''_{g + g' + i}) \leq D(\Sigma_{g+i}) + D(\Sigma_{g'})$.
Because $D(\Sigma_{g+i}) = D(\Sigma_g)$ for all $i$, we have  
$A(M \# M') \leq A(M) + A(M')$.
\end{proof}

The converse statement is an open question: Is it
necessarily true that $A(M \# M') \geq A(M) + A(M')$?

\begin{Lem}
\label{dehnfillem}
Let $M$ be a manifold with boundary and let $M'$ be the result of filling
one or more torus boundary components with solid tori. Then $A(M') \leq A(M)$.
\end{Lem}

\begin{proof}
Let $\Sigma$ be a Heegaard surface for $M$ and let $\Sigma'$ be the image
of $\Sigma$ in the induced map $M \rightarrow M'$.  This $\Sigma'$ is a
Heegaard surface for $M'$.  Let $v_1,\dots,v_n$ be a path in $C^*(\Sigma)$
from $H_1$ to $H_2$.  The map $\Sigma \rightarrow \Sigma'$ suggests an
isomorphism $C^*(\Sigma) \rightarrow C^*(\Sigma')$.
Let $v'_1,\dots,v'_n$ be the images of $v_1,\dots,v_n$.

The vertex $v_n$ defines $H'_2$ because $H_2$ is a handlebody so $H'_2$
is the image of $H_2$ in the induced map.  However, $v_1$ may not define
$H_1$ because for each torus boundary component that is filled, there is a
loop in $\Sigma$ parallel to $\partial_- H_1$ which may not bound a disk
in $H'_1$.

Let $T$ be a torus boundary component which is filled so that the boundary 
of a meridian disk maps to a loop $\alpha \subset T$.
If we cut $H_1$ along the maximal collection disks defined by the pants
decomposition, there will be a loop $l_1$ on the component containing $T$
and a scar bound by a loop $l_2$.  The loop $\alpha$ (which sits in 
$\partial_- H_1$) can be projected into $\partial_+ H_1$ so that its image
sits in the torus with $l_1$ and is disjoint from $l_2$.  Thus it takes one 
move to replace $l_1$ with the image of $\alpha$.  

The image of $\alpha$ bounds a disk in $H'_1$.  For each torus component
which is filled, it takes at most one move to replace a loop in the pants 
decomposition for $H_1$ with a loop bounding a disk in $H'_1$.  The final
product is a pants decomposition containing a maximal collection of disks
for $H'_1$.  The sum $b'$ of the genera of the boundary components goes down
by one for every Dehn filling.

Since we have $A_g(\Sigma') = D(\Sigma'_g) - g + b'$, we know that 
$A_g(\Sigma') \leq A_g(\Sigma)$ when $g$ is the genus of $\Sigma$.  The
same proof works for every stabilization of $\Sigma$ so in the limit
we have $A(M') \leq A(M)$.
\end{proof}

The equivalent statement is not true for $A^P(M)$ because 
there is no control over the number of times the image of the loop $\alpha$
intersects the loop $l_1$.  

Although the Heegaard complexity cannot increase under Dehn filling, 
it can drop by an arbitrary amount.  In particular, if $M$ is the complement
of a knot in $S^3$, then there is a Dehn filling which produces $M' = S^3$, so 
$A(M') = 0$, regardless of the Heegaard
complexity of $M$.

\begin{Lem}
\label{irredbndlem1}
Let $(\Sigma, H_1, H_2)$ be an irreducible Heegaard splitting of a closed
manifold $M$ and let $g$ be the genus of $\Sigma$.  Then 
$D(\Sigma) \geq 3g - 3$.
\end{Lem}

\begin{proof}
Let $v_0,\dots,v_n$ be a path in $C^*(\Sigma)$ such that $v_1$ defines $H_1$,
$v_n$ defines $H_2$ and $n = D(\Sigma)$.  Both $H_1$ and $H_2$ are 
handlebodies so if $\sigma_{v_1}$ and 
$\sigma_{v_n}$ share a vertex $u \in C(\Sigma)$ then $u$ corresponds to a 
loop in $\Sigma$ that bounds a disk in $H_1$ and a disk in $H_2$.  

Since 
$\Sigma$ is irreducible, there is no such loop in $\Sigma$
so $\sigma_{v_0}$ and $\sigma_{v_n}$ cannot share a vertex.  The simplex
$\sigma_{v_0}$ has $3g - 3$ vertices and consecutive simplices
$\sigma_{v_i}, \sigma_{v_{i+1}}$ share all but one vertex so there must
be at least $3g-2$ in the sequence.  Thus it must be that $n \geq 3g-2$
and $D(\Sigma) \geq 3g - 3$.
\end{proof}

\begin{Lem}
\label{irredbndlem2}
Let $(\Sigma, H_1, H_2)$ be a strongly irreducible Heegaard splitting of a 
closed manifold $M$ and let $g$ be the genus of $\Sigma$.  Then 
$D(\Sigma) \geq 6g - 7$.
\end{Lem}

We will sketch the proof, since the result is not vital to the rest of the
paper.  Because $\Sigma$ is irreducible, each of the loops must be moved at 
least once.  Let $l$ be the last loop that's moved.  A loop $l'$ which is moved
before $l$ is replaced by a loop $l''$ disjoint from $l$.  Since $l$ bounds a 
disk in $H_1$, $l''$ cannot bound a disk in $H_2$ (since $\Sigma$ is strongly
irreducible) so $l''$ must be moved later on.  Every loop other than $l$
must therefore be moved at least twice (once before $l$ and once after)
so $D(\Sigma) \geq 6g - 7$.

Unfortunately, once the Heegaard surface $\Sigma$ is replaced with a
stabilization $\Sigma'$, the distance may drop by an arbitrary amount.
For example, Kobayashi~\cite{kob:poly} has constructed a manifold, based
on work by Casson and Gordon, with a sequence of strongly irreducible
Heegaard splittings of arbitrarily high genus.  Sedgwick~\cite{sedg:stab} 
later showed that the result of stabilizing any of these once is also a 
stabilization of all the lower genus splittings.

By \fullref{irredbndlem2}, the unstabilized Heegaard splittings have
arbitrarily high distance, but
by \fullref{uprbndlem}, the distances of the stabilizations are
bounded, so the distance must fall by an arbitrarily large amount after
stabilization.

There is still a relationship between the Heegaard genus of the manifold
and the Heegaard complexity, but it is not as strong.

\begin{Lem}
\label{genbndlem}
If $M$ is irreducible and $\partial M = \emptyset$ then the Heegaard genus of
$M$ is less than or equal to $\frac{A(M)+2}{2}$.
\end{Lem}

\begin{proof}
Let $\Sigma$ be a Heegaard splitting of $M$ such that $A_g(\Sigma) = A(M)$,
where $g$ is the genus of $\Sigma$.  Let $v_1,\dots,v_n$ be a path in
$C^*(\Sigma)$ of length $D(\Sigma)$.  Let $\mathbf{L}$ be the loops in $\Sigma$
defined by the vertices of $\sigma_{v_1}$.  By choosing $g$ large enough,
we can assume there is at least one loop $l_i$ which also bounds a disk in
$H_2$.

Let $\mathbf{L'}$ be the collection of loops which are never moved.  In other
words if $l_j \in \mathbf{L'}$ then for every vertex $v_i$ in the sequence,
$l_j$ corresponds to a vertex of $\sigma_{v_i}$.  In particular, $l_j$
is a vertex of $\sigma_{v_1}$ and of $\sigma_{v_n}$.  Thus $l_j$ bounds
a disk $D^1_j \subset H_1$ and a disk $D^2_j \subset H_2$ (since $H_1$ and 
$H_2$ are handlebodies).  We can choose the disks corresponding to the loops 
in $\mathbf{L}$ so that the collection is pairwise disjoint.  Thus each $l_j$
suggests a disjoint, embedded 2--sphere $S_j$.

Because $M$ is irreducible, each sphere $S_j$ bounds a ball $B_j \subset M$.
Let $\Sigma' = \Sigma \setminus \bigcup B_j$.  This is a punctured surface.
Disks can be glued into the punctures to make $\Sigma'$ a Heegaard surface
for $M$, so the genus of $\Sigma'$ is at least the Heegaard genus of $M$.
Let $h$ be the genus of $\Sigma'$.  Some subset $\mathbf{L''}$ of the loops 
$\mathbf{L}$ form a pair-of-pants decomposition for $\Sigma'$.  This is a
punctured surface of genus $h$ (there may be more than one puncture) so
there are at least $3h - 2$ loops in $\mathbf{L''}$ and none of these
loops are in $\mathbf{L'}$.

In the surface $\Sigma \setminus \Sigma'$, there are at least $g - h$
non-separating loops.  Since $M$ is irreducible, none of these loops can
be in $\mathbf{L'}$ so there are at least $g - h$ more loops that are
moved.  There are at least $2h - 2 + g$ loops in 
$\mathbf{L} \setminus \mathbf{L'}$ so $D(\Sigma) \geq 2h - 2 + g$
and $A_g(\Sigma) \geq 2h - 2$.  This is true of every stabilization of
$\Sigma$ so $h \leq \frac{A(M)+2}{2}$ and the Heegaard genus of $M$ is
at most $h$. 
\end{proof}

\begin{Coro}
For every positive integer $N$, there is a manifold $M$ with
$A(M) > N$.
\end{Coro}

\section{Comparing manifolds}
\label{compmfldsect}

Let $M$ and $M'$ be compact, connected, orientable 3--manifolds such that
$\partial M = \partial M'$.  (Both boundaries may be empty.) Let
$(\Sigma, H_1, H_2)$ and $(\Sigma', H_1', H_2')$ be Heegaard splittings of
$M$, $M'$ respectively such that $H_2$ and $H_2'$ are handlebodies of
the same genus.

Let $\phi \co H_1 \rightarrow H_1'$ be any homeomorphism.  Such a map exists
because $\partial_- H_1 =\partial M= \partial M' =\partial_- H_1'$.
This induces a homeomorphism from $\Sigma = \partial_+ H_1$ to
$\Sigma' = \partial_+ H_1'$ and this homeomorphism suggests an isomorphism
$\hat \phi \co C^*(\Sigma) \rightarrow C^*(\Sigma')$ or
$\hat \phi \co C^P(\Sigma) \rightarrow C^P(\Sigma')$.

\begin{Def}  
The \emph{dual distance} between the two Heegaard splittings is 
$D(\Sigma, \Sigma')=\{D(v, v') : v \text{ defines } H_2,\, \hat\phi(v')
\text{ defines } H_2'\}$.  The \emph{pants distance} is
$D^P(\Sigma, \Sigma')=\min \{D^P(v, v') : v \text{ defines } H_2,\,\hat\phi(v')
\text{ defines } H'_2\}$.
Both minima are taken over all homeomorphisms $\phi \co H_1 \rightarrow H_1'$.
\end{Def}

In other words, identify $H_1$ and $H_1'$ and consider two pants
decompositions of $\partial_+ H_1$ such that one is a Heegaard diagram for $M$
and the other is a Heegaard diagram for $M'$.  The value of
$D(\Sigma, \Sigma')$ is the smallest possible dual distance between any two
such pants decompositions.

As in \fullref{stabsect}, consider a sequence $\Sigma_g$ of genus--$g$
stabilizations of $\Sigma$ and a sequence $\Sigma'_g$ of stabilizations
of $\Sigma'$.

\begin{Lem}
$D(\Sigma_g, \Sigma'_g) \leq 2 D(\Sigma, \Sigma')$.
\end{Lem}

The proof is almost identical to the proof of \fullref{uprbndlem},
with the exception that there are not $(h - g)$ loops left at the end
which need to be moved.  The proof is left to the reader, as is the proof
of the following lemma.

\begin{Lem}
For sufficiently large $g$, $D(\Sigma_h,\Sigma_h') \leq D(\Sigma_g,\Sigma'_g)$
whenever $h \geq g$.
\end{Lem}

\begin{Thm}
The sequence $D(\Sigma_g, \Sigma_g')$ converges and is independent of the
choices of $\Sigma$ and $\Sigma'$.
\end{Thm}

\begin{proof}
The distance $D(\Sigma_g, \Sigma_g')$ is non-negative and the sequence is
non-increasing for sufficiently large $g$, so it converges.  Showing that
the limit is unique is slightly more tricky.

Let $(\Sigma^1, H_1^1, H_2^1)$ and $(\Sigma^2, H_1^2, H_2^2)$ be Heegaard 
splittings of $M$.  We know that there is a stabilization 
$(\Sigma^3, H_1^3, H_2^3)$ of $\Sigma^1$ and a stabilization
$(\Sigma^4, H_1^4, H_2^4)$ of $\Sigma^2$ such that $\Sigma^3$ is isotopic
to $\Sigma^4$.  If $M$ has non-empty boundary then this isotopy must send
$H_1^3$ to $H_1^4$ since $H_2^4$ is a handlebody but $H_1^3$ is not.

However, 
if $M$ is closed, the isotopy could send $H_1^3$ to $H_2^4$.  Since the
definition of $D(\Sigma, \Sigma')$ distinguishes between $H_1$ and $H_2$,
the distances $D(\Sigma^3, \Sigma')$ and $D(\Sigma^4, \Sigma')$ may not
be equal.

Implicit in most proofs of the stabilization theorem it is also proven
that the stabilizations can be chosen so that an isotopy sends 
$\Sigma^3$ to $\Sigma^4$ and sends $H^3_1$ to $H^4_1$.  (See, for example,
Rubinstein and Scharlemann's proof~\cite{rub:compar}.)  In this case, it
must be the case that $D(\Sigma^3, \Sigma') = D(\Sigma^4, \Sigma')$.
So the limit of the distances is independent of the choice of $\Sigma$.
A similar consideration for $M'$ shows that the distance is independent of 
the choice of $\Sigma'$, and therefore the limit depends only on $M$ and $M'$.
\end{proof}

\begin{Def}
We define the \emph{Heegaard distance} $D(M, M')$ to be the limit of
the sequence $D(\Sigma_g, \Sigma_g')$ for any Heegaard surfaces $\Sigma$
and $\Sigma'$ of $M$ and $M'$.
\end{Def}

It is immediate that $D(M, M') = 0$ if and only if $M = M'$.  When 
$M \neq M'$, there is a very simple description of $D(M, M')$ as  follows:  
Let $K \subset M$ be a link.  We will say that $M$ and $M'$ are connected by
$K$ if $M'$ is the result of some Dehn surgery on $K \subset M$.
Let $c(M, M')$ be the smallest integer $c$ such that $M$ and $M'$ are
connected by an $c$--component link.

\begin{Thm}
\label{dehnthm}
$D(M, M') = c(M, M')$.
\end{Thm}

\begin{proof}
We will first show that $D(M, M') \leq c(M, M')$.
Let $K \subset M$ be a $c$--component link in $M$ such that some Dehn surgery
on $K$ yields $M'$.  Let $(\Sigma, H_1, H_2)$ be a Heegaard splitting of
$M$ such that there is a trivalent spine $G$ of $H_2$ in which each component 
of $K$ appears as an edge of $G$ with both ends on the same vertex.

Let $l_1,\dots,l_m$ be loops on $\Sigma$ defined by the meridian disks dual
to the edges of $G$ and assume
$l_1,\dots,l_n$ correspond to the edges of $G$ defined by the components of
$K$.  Since $G$ is trivalent, the loops suggest a pants decomposition of
$\Sigma$ and the corresponding vertex $v_0$ of $C^*(\Sigma)$ defines $H_2$.

The loop $l_1$ defines a meridian curve on a component of $K$.  Let $l_1'$
be the meridian defined by the Dehn surgery on $K$.  This loop is disjoint
from the loops $l_2,\dots,l_m$ because of the way we chose $G$.  Thus
the vertex $v_1$ defined by $l_1',l_2,\dots,l_m$ is connected to $v_0$ by
an edge in $C^*(\Sigma)$.

Continuing in this way for the loops $l_2,\dots,l_n$, we can construct a
path $v_0,\dots,v_n$ where $v_n$ is defined by the loops
$l_1',\dots,l_n',l_{n+1},\dots,l_m$.  By the construction, this vertex will
define a handlebody in a Heegaard splitting $(\Sigma',H_1', H_2')$ of $M'$.
For any stabilization of $(\Sigma, H_1, H_2)$, we can construct a graph with
the same properties as $G$.  Thus $D(\Sigma_g, \Sigma'_g) \leq c$ for every
$g$ and $D(M, M') \leq c$.

We now see a connection between Heegaard distance and Dehn surgery.  To
prove that $D(M,M') \geq c$, we need the following Lemma:

\begin{Lem}
\label{oneknotlem}
If $D(M, M') = 1$ then $M$ and $M'$ are connected by a knot (a one-component
link).
\end{Lem}

\begin{proof}
Let $(\Sigma, H_1, H_2)$ and $(\Sigma',H_1',H_2')$ be Heegaard splittings for
$M$, $M'$ respectively such that $D(\Sigma, \Sigma') = 1$.  This means there
is a pants decomposition $l_1,l_2,\dots,l_m$ for $\partial H_1$ giving a
Heegaard diagram for $M$ such that replacing $l_1$ with $l'_1$ creates a 
pants decomposition $l_1',l_2,\dots,l_m$ giving a Heegaard diagram for $M'$.

The loop $l_1$ bounds a disk in $H_2$ and the loop $l_1'$ bounds a disk in
$H_2'$.  As in \fullref{defsect}, $l_1$ may sit in a four-punctured
sphere or a once-punctured torus, and $l_1'$ will lie in the same type
of component.  If $l_1$ sits in a four-punctured sphere, then the punctured
sphere and the disks defined by the punctures bound a ball in $H_2$.
Since $l_1'$ is on the boundary of this ball, it bounds a disk in $H_2$.
Thus $l_1',l_2,\dots,l_m$ is a Heegaard diagram for $M$ as well as $M'$ so
the two manifolds are homeomorphic.  Since $D(M, M') \neq 0$, the manifolds 
are distinct, and $l_1$ must sit in a punctured torus.

The punctured torus and the disk defined by the puncture bound a solid
torus in $M$ and a solid torus in $M'$.  The remainder of the Heegaard
diagrams are identical so the complements of the solid tori in $M$ and
$M'$ are homeomorphic.  We get $M$ and $M'$ by gluing solid tori to
the boundaries of the complement, so Dehn surgery on the solid torus in
$M$ will yield $M'$.
\end{proof}

We can now finish the proof of \fullref{dehnthm}

Let $v_1,\dots,v_n$ be a path in $C^*(\Sigma)$ from $H_2$ to $H_2'$.  Each
vertex $v_i$ defines a Heegaard splitting of a manifold $M_i$.  Because the
Heegaard distance between consecutive manifolds is at most 1, there is a 
sequence of knots such that Dehn surgery on each knot yields the next 
manifold.  By keeping track of the images of these knots in $M$, we find 
a link with at most as many components as the distance of the path.  Thus
$c(M,M') \leq D(M,M')$ and the proof is complete.
\end{proof}

As one might expect, there is an analogous Theorem for $D^P(M,M')$.
Let $c^1(M, M')$ be the minimum number of components of a link $K \subset M$
such that $M'$ is the result of a surgery in which a meridian of each 
component is replaced by a loop which intersects the meridian once.  

\begin{Thm}
$D^P(M, M') = c^1(M, M')$.
\end{Thm}

The proof is almost identical to the proof of \fullref{dehnthm} and will
be left to the reader.

\section{Questions and speculations}
\label{questsect}

\begin{Qu}
Is $A(M)$ or $A^P(M)$ related to (quasi-equal to?) a manifold invariant which 
is already known?
\end{Qu}

The most tempting possibility is that for hyperbolic manifolds, the
complexity could be related to volume.  Juan Souto~\cite{souto:volume} has 
announced joint
work with Jeff Brock showing that in certain cases, the pants distance of a
specific Heegaard splitting is quasi-isometric to its hyperbolic
volume.  Brock~\cite{Brock:pants} has previously proven results relating 
distance in the pants complex to volumes of convex cores.

Another possibility is that the pants complexity could be
related to the number of tetrahedra needed for a triangulation of $M$, or
to triangulate a 4--manifold bounded by $M$.  By Hatcher and Thurston,
an edge path in the pants complex suggests a smooth path in 
$C^{\infty}(\Sigma)$.  
Applying these functions to the level surfaces of
a sweep-out of the Heegaard splitting suggests a stable function from the
manifold to $\mathbb{R}^2$.  Costantino and Thurston~\cite{thr:shadows} have 
used 2--dimensional stable functions to relate the number of tetrahedra in a
triangulation of $M$ to the number of tetrahedra in a 4--manifold which it 
bounds.  The link between the complexity of the stable function and
triangulations may be applicable to the stable function induced by a
Heegaard splitting.

\begin{Qu}
Is there any way to calculate the complexity?
\end{Qu}

This is already known to be a very difficult problem.  There is currently
a lot of work being done to calculate the distance between two points
in the pants complex or the Hempel distance of a Heegaard splitting.
Calculating the Heegaard complexity should be even more difficult because
it requires calculating the distance for an infinite number of splittings.

\begin{Qu}
Is either complexity additive under connect summing?
\end{Qu}

This is a deceptively simple-sounding problem.  We have seen that both 
complexities are sub-additive.  Unlike the Hempel distance, which is
zero when $M$ is reducible, the dual distance or pants distance may
ignore the reducing disks when finding a shortest path.

\begin{Qu}
What conditions will guarantee that $A^P(M)$ will not increase after
Dehn filling a torus boundary component?  What conditions will guarantee
that $A(M)$ will not decrease?
\end{Qu}

It was pointed out earlier that there is no analogy to \fullref{dehnfillem} for
$A^P(M)$. Given a manifold with boundary, a Heegaard splitting and a minimal 
distance path in $C^P(\Sigma)$, the initial pants decomposition of $\Sigma$ 
will define loops in $\partial M$.  Dehn fillings which take into account these
loops will guarantee that $D^P(\Sigma)$ does not increase.  The question is
whether or not there is a way to predict these loops from the topology of
the manifold, ie without calculating geodesic paths for an infinite
sequence of Heegaard splittings.

\begin{Qu}
How does the complexity behave under finite covers?
\end{Qu}

Lackenby~\cite{Lack:asym} has shown that the asymptotic behavior of the 
Heegaard genus of finite covers is related to Thurston's virtually Haken 
conjecture and
virtually fibered conjecture.  Essentially, if the Heegaard genera
of finite covers of a given manifold are bounded by a nice enough function
of the degree of the covers then one of the covers must be Haken or
fibered.

Given a finite cover $M'$ of $M$ and a Heegaard splitting, there is an 
induced Heegaard splitting of $M'$.  By  \fullref{genbndlem}, rather than 
having to find an
alternative splitting of $M$ with lower genus, it would only
be necessary to show that the pants distance is bounded by a nice enough
function of the degree of the cover.

\begin{Qu}
If $M$ and $N$ are irreducible 3--manifolds and there is a degree-one map from 
$M$ to $N$, does this imply that $A(N) \leq A(M)$?
\end{Qu}

Waldhausen~\cite{wald:degone} showed that given a degree-one map 
$f \co M \rightarrow N$ and a Heegaard splitting $(\Sigma, H_1, H_2)$ of $N$, 
one can construct a Heegaard splitting $(\Sigma', H'_1, H'_2)$ of $M$ such that
$f$ maps $H'_1$ and $H'_2$ onto $H_1$ and $H_2$, respectively, by a simple type
of degree-one map.  Given a path in $C^*(\Sigma)$, it may be possible to 
construct a path in $C^*(\Sigma')$ which is, in some sense, induced by the map 
$f$.

A positive answer to this question would imply the Poincare conjecture because 
a homotopy equivalence is a degree-one map.  So, if $M$ is $S^3$ then 
$A(M) = 0$ so $A(N) = 0$ and $N$ is $S^3$.  If we remove the assumption that 
$M$ is irreducible then the answer is no, since there is a degree-one map from 
$S^1 \times S^2$ to $S^3$.

\bibliographystyle{gtart}
\bibliography{link}

\end{document}